\documentclass[12pt, twoside, leqno]{article}
\usepackage{amsmath,amsthm}
\usepackage{amssymb,latexsym}
\usepackage{enumerate}

\newtheorem{thm}{Theorem}[section]
\newtheorem{prop}[thm]{Proposition}
\newtheorem{cor}[thm]{Corollary}
\newtheorem{lem}[thm]{Lemma}

\newtheorem{conj}[thm]{Conjecture}

\theoremstyle{definition}

\newtheorem{rem}[thm]{Remark}

\numberwithin{equation}{section}

\usepackage{mathrsfs}
\usepackage{CJK}
\usepackage{amsmath, amscd}
\usepackage{fancyhdr}
\usepackage{enumerate}

\begin{document}

\baselineskip=17pt

\title{Simple derivations and their images}
\author{
Ruiyan Sun\\
MOE-LCSM,\\ School of Mathematics and Statistics,\\
Hunan Normal University, Changsha 410081, China \\
\emph{E-mail:} sunruiyan19@163.com \\
Dan Yan\footnote{ The author is supported by Scientific Research Fund of Hunan Provincial Education Department (Grant No. 21A0056), the NSF of China (Grant No. 11871241; 11601146) and the Construct Program of the Key Discipline in Hunan Province.}\\
MOE-LCSM,\\ School of Mathematics and Statistics,\\
Hunan Normal University, Changsha 410081, China \\
\emph{E-mail:} yan-dan-hi@163.com \\
}
\date{}
\maketitle

\renewcommand{\thefootnote}{}

\renewcommand{\thefootnote}{\arabic{footnote}}
\setcounter{footnote}{0}

%%%%%%%%
\begin{abstract}
In the paper, we prove that the derivation $D=y\partial_x+(a_2(x)y^2+a_1(x)y+a_0(x))\partial_y$ of $K[x,y]$ with $a_2(x),a_1(x),a_0(x)\in K[x]$ is simple iff the following conditions hold:
$(1)$ $a_0(x)\in K^*$,
$(2)$ $\deg a_1(x)\geq1$ or $\deg a_2(x)\geq1$,
$(3)$ there exist no $l\in K^*$ such that $a_2(x)=la_1(x)-l^2a_0(x)$.
In addition, we prove that the image of the derivation $D=\partial _x+{\sum_{i=1}^n \gamma_i(x) y_i^{k_i}}{\partial _i}$ is a Mathieu-Zhao space iff $D$ is locally finite. Moreover, we prove that the image of the derivation $D={\sum_{i=1}^n \gamma_i y_i^{k_i}}{\partial _i}$ of $K[y_1,\ldots,y_n]$ is a Mathieu-Zhao space iff $k_i\leq 1$ for all $1\leq i\leq n$, $n\geq 2$.
\end{abstract}
{\bf Keywords.} Simple derivations, Darboux polynomials, Mathieu-Zhao spaces\\
{\bf MSC(2020).} 13N15; 14R10; 13P05. \vskip 2.5mm

\baselineskip=17pt

\section{Introduction}

Throughout this paper, we will write $K$ for any field with
characteristic zero, $K^*$ for the set of all elements satisfy that $a\in K$ and $a\neq0$, $\mathbb{N}$ for the set of all natural numbers including zero and $R:=K[x,y_1,\dots,y_n]$ for the
polynomial algebra over $K$ in $n+1$ indeterminates $x,y_1,\ldots,y_n$.
$\partial_x,~\partial_i$ will denote the derivations $\frac{\partial}{\partial x}$, $\frac{\partial}{\partial y_i}$ of $R$ for all $1\leq i\leq n$, respectively. For an element $f$ of $K[z]$, we often use $f'$ instead of $f_z$, where $z$ is an indeterminate over $K$ and $z\in\{x,y_1,\ldots,y_n\}$.

A $K$-derivation $D:R\rightarrow R$ of $R$ is a $K$-linear map such that
$$D(ab)=D(a)b+aD(b)$$
for any $a, b\in R$ and $D(c)=0$ for any $c\in K$. An ideal $I$ of $R$ is called $D$-stable if $D(I)\subset I$. $R$ is called
$D$-simple if it has no proper nonzero $D$-stable ideal. The $K$-derivation $D$ is called simple if $R$ has no $D$-stable
ideal other than $0$ and $R$. For some examples of simple derivations, see \cite{2}, \cite{4}, \cite{1}, \cite{7}, \cite{15}.

A polynomial $F\in K[x,y]$ is said to be a Darboux polynomial of $D$ if $F\notin K$ and $D(F)=\Lambda F$ for some $\Lambda\in K[x,y]$. We define $\deg 0=-\infty$.

Let $D$ be any $K$-derivation of $R$. Then $D$ is said to be locally finite if for every $a\in R$ the $K$-vector space spanned by the elements $D^ia (i\geq 1)$ is finite dimensional.

The Mathieu-Zhao spaces were introduced by Zhao in \cite{18}. We give the definition here for the polynomial rings. A $K$-subspace $M$ of $R$ is said to be a Mathieu-Zhao space if $a\in R$ such that $a^m\in M$ for all large $m$, then for every $b\in R$, we have $ba^m\in M$ for all large $m$.

Simple derivations are useful to construct simple noncommutative rings which are not fields(\cite{14}). If $n\geq 1$, then only some examples of simple derivations of $R$ are known.
 A. Nowicki in \cite{8} or A. Shamsuddin in \cite{3} gave an algorithm to decide whether a $K$-derivation of $K[x,y]$ of the form $D=\partial_x+(a(x)y+b(x))\partial_y$, where $a(x), b(x)\in K[x]$, is simple or not. In \cite{5}, Y. Lequain has characterized the property ``$D$ is simple of $R$'' in terms of certain properties of $D$ that one can effectively check for whether it is satisfied or not. In \cite{17}, S. Kour has shown that $D=y^r\partial_x+(xy^s+g)\partial_y$, where $0\leq r<s$ are integers, $g\in K[y]$, is a simple derivation of $K[x,y]$. In \cite{15}, the second author proved that $D=y\partial_x+(a_2y^2+a_1(x)y+a_0(x))\partial_y$ with $a_2\in K$, $a_1(x), a_0(x)\in K[x]$ is simple iff $a_0(x)\in K^*$ and $\deg a_1(x)\geq 1$.

In our paper, we give a necessary and sufficient condition for a derivation to be simple and study the images of some derivations. In section 2, we prove that $D=y\partial_x+(a_2(x)y^2+a_0(x))\partial_y$ with $a_2(x),a_0(x)\in K[x]$ is simple iff $a_0(x)\in K^*,\deg a_2(x)\geq 1$. In section 3, we prove some lemmas which we need in section 4. In section 4, we prove that $D=y\partial_x+(a_2(x)y^2+a_1(x)y+a_0(x))\partial_y$ with $a_2(x),a_1(x),a_0(x)\in K[x]$ is simple iff $(1)$ $a_0(x)\in K^*$,
$(2)$ $\deg a_1(x)\geq 1$ or $\deg a_2(x)\geq 1$,
$(3)$ there exist no $l\in K^*$ such that $a_2(x)=la_1(x)-l^2a_0(x)$.
 Finally we prove that the image of the derivation $D=\partial _x+{\sum_{i=1}^n \gamma_i(x) y_i^{k_i}}{\partial _i}$, where $k_i\geq 1, \gamma_i(x)\in K[x]$ for all $1\leq i\leq n$, is a Mathieu-Zhao space iff $D$ is locally finite. In addition, we show that the image of the derivation $D={\sum_{i=1}^n \gamma_i y_i^{k_i}}{\partial _i}$, where $n\geq 2$ and $\gamma_i\in K ^*$ for all $1\leq i\leq n$, is a Mathieu-Zhao space iff $k_i\leq1$ for all $1\leq i\leq n$. Base on the conclusions of section 2 to section 4, we give a conjecture in section 6.

\section{The case $D=y\partial_x+(a_2(x)y^2+a_0(x))\partial_y$}

\begin{thm} \label{thm2.1}
Let $D=y\partial_x+(a_2(x)y^2+a_0(x))\partial_y$ be a derivation of $K[x,y]$ with $a_2(x),a_0(x)\in K[x]$. Then $D$ is simple if and only if $a_0(x)\in K^*$ and $\deg a_2(x)\geq1$.
\end{thm}
\begin{proof}
$``\Rightarrow"$  If $a_0(x)\notin K^*$, then $I=(y,a_0(x))$ is a $D$-stable ideal. Thus, $D$ is not simple. Hence we have $a_0(x)\in K^*$. If $\deg a_2(x)\leq 0$, then $a_2(x)\in K$. If $a_2(x)\in K^*$, then $I=(y^2+\frac{a_0(x)}{a_2(x)})$ is a $D$-stable ideal. Thus, $D$ is not simple. If $a_2(x)=0$, then $I=(\frac{1}{2}y^2-xa_0(x))$ is a $D$-stable ideal. Thus, $D$ is not simple. Hence we have $\deg a_2(x)\geq1$.

$``\Leftarrow"$ It follows from Proposition 2.1 in \cite{17} that it suffices to prove that $(y,a_2(x)y^2+a_0(x))=(1)$ and $D$ has no Darboux polynomial. Clearly, $(y,a_2(x)y^2+a_0(x))=(a_0(x))=(1)$. Suppose that $F(x,y)\in K[x,y]$ is a Darboux polynomial of $D$ and $F(x,y)={\sum_{i=0}^n c_i(x) y^i}\notin K$ with $c_n(x)\neq 0,c_i(x)\in K[x]$ for $0\leq i \leq n$. Then we have $D(F(x,y))=\Lambda(x,y)F(x,y)$ for some $\Lambda(x,y)\in K[x,y]$.

If $n=0$, then $F(x,y)=c_0(x)$, $yc'_0(x)=\Lambda(x,y)c_0(x)$. Thus, we have $c_0(x)\in K$, which is a contradiction. Hence we have $n\geq1$. It follows from the equation $D(F(x,y))=\Lambda(x,y)F(x,y)$ that
\begin{equation}
y(\sum_{i=0}^n{c'_i(x)}y^i)+(a_2(x)y^2+a_0(x))(\sum_{i=1}^n{i}c_i(x)y^{i-1})=\Lambda(x,y)\sum_{i=0}^nc_i(x)y^i. \label{eq2.1}
\end{equation}

Comparing the highest degree of $y$ of equation \eqref{eq2.1}, we have $\deg_y\Lambda(x,y)\leq1$. Let $\Lambda(x,y)=d_1(x)y+d_0(x)$ with $d_1(x),d_0(x)\in K[x]$. Then equation \eqref{eq2.1} has the following form:
\begin{equation}
y(\sum_{i=0}^n{c'_i(x)}y^i)+(a_2(x)y^2+a_0(x))(\sum_{i=1}^n{i}c_i(x)y^{i-1})=(d_1(x)y+d_0(x))\sum_{i=0}^nc_i(x)y^i. \label{eq2.2}
\end{equation}

We view the polynomials as in $K[x][y]$ with coefficients in $K[x]$ when we comparing the coefficients of monomials on $y$. Comparing the coefficients of $y^{n+1}$ of equation \eqref{eq2.2}, we have
\begin{equation}c'_n(x)+na_2(x)c_n(x)=d_1(x)c_n(x).\label{eq2.10}
\end{equation}

Thus, we have $c_n(x)\in K^*$ and $d_1(x)=na_2(x)$ by comparing the highest degree of $x$ of the above equation. If $n\geq2$, comparing the coefficients of monomials on $y$ of equation (2.2), then we have the following equations:
\begin{eqnarray}
c'_{n-1}(x)&=&a_2(x)c_{n-1}(x)+d_0(x)c_n(x),  \label{eq2.3} \\
(i+1)a_0(x)c_{i+1}(x)&=&(n-i+1)a_2(x)c_{i-1}(x)+d_0(x)c_i(x)-c'_{i-1}(x), \label{eq2.4} \\
a_0(x)c_1(x)&=&d_0(x)c_0(x) \label{eq2.5}
\end{eqnarray}
for $1\leq i\leq{n-1}$.

If $n=1$, comparing the coefficients of monomials on $y$ of equation \eqref{eq2.2}, then we have equations \eqref{eq2.3},\eqref{eq2.5}.

Claim 1: $d_0(x)\neq0$.

Suppose that $d_0(x)=0$. It follows from equation \eqref{eq2.3} that $a_2(x)c_{n-1}(x)-c'_{n-1}(x)=0$. Thus, we have $c_{n-1}(x)=0$ by comparing the degree of $x$ of equation \eqref{eq2.3}. If $n=1$, then we have $c_0(x)=0$. It follows from equation \eqref{eq2.5} that $a_0(x)c_1(x)=0$, which is a contradiction. If $n\geq 2$, then it follows from equation \eqref{eq2.4}$(i=n-1)$ that $na_0(x)c_n(x)=2a_2(x)c_{n-2}(x)-c'_{n-2}(x)$. Since $\deg a_2(x)\geq1$ and $na_0(x)c_n(x)\in K^*$, we have a contradiction by comparing the degree of $x$ of equation \eqref{eq2.4} $(i=n-1)$. Hence we have $d_0(x)\neq0$.

It follows from equation \eqref{eq2.3} that $c_{n-1}(x)\neq0$. Otherwise, we have $d_0(x)c_n(x)\\=0$, which is a contradiction. Thus, it follows from equation \eqref{eq2.3} that
\begin{equation}
\deg d_0(x)=\deg a_2(x)+\deg c_{n-1}(x)\geq1.
\label{eq2.6}
\end{equation}

 If $n\geq2$, then it follows from equation \eqref{eq2.4}$(i=n-1)$ that $c_{n-2}(x)\neq0$. Otherwise, $na_0(x)c_n(x)=c_{n-1}(x)d_0(x)$. Then we have a contradiction by comparing the degree of $x$ of equation \eqref{eq2.4}$(i=n-1)$. Thus, it follows from equation \eqref{eq2.4}$(i=n-1)$ that $\deg c_{n-2}(x)+\deg a_2(x)=\deg c_{n-1}(x)+\deg \deg d_0(x)$. It follows from equation \eqref{eq2.6} that $\deg c_{n-2}(x)=2\deg c_{n-1}(x)$.

Claim 2: $c_{n-i}(x)\neq0$ and $\deg c_{n-i}(x)=i\deg c_{n-1}(x)$ for all $1\leq i\leq n$.

Suppose that $c_{n-l}(x)\neq0$ and $\deg c_{n-l}(x)=l\deg c_{n-1}(x)$ for all $1\leq l\leq k-1,k\in \{2,3,\ldots,n\}$. Next we show that $c_{n-k}(x)\neq0$ and $\deg c_{n-k}(x)=k\deg c_{n-1}(x)$.

If $c_{n-k}(x)=0$, then it follows from equation \eqref{eq2.4} $(i=n-k+1)$ that
\begin{equation}
(n-k+2)a_0(x)c_{n-k+2}(x)=d_0(x)c_{n-k+1}(x).\label{eq2.7}
\end{equation}

It follows from equation \eqref{eq2.7} that
\begin{equation}
\deg c_{n-k+2}(x)=\deg d_0(x)+\deg c_{n-k+1}(x).\label{eq2.8}
\end{equation}

By induction hypothesis, we have $\deg c_{n-k+2}(x)=(k-2)\deg c_{n-1}(x)$ and $\deg c_{n-k+1}(x)=(k-1)\deg c_{n-1}(x)$. Then equation \eqref{eq2.8} has the following form:
$\deg d_0(x)+\deg c_{n-1}(x)=0$, which contradicts equation \eqref{eq2.6}. Hence we have $c_{n-k}(x)\neq0$. It follows from equation \eqref{eq2.4}$(i=n-k+1)$ that
\begin{equation}
 \deg a_2(x)+\deg c_{n-k}(x)=\deg d_0(x)+\deg c_{n-k+1}(x).\label{eq2.9}
\end{equation}

Then we have $\deg c_{n-k}(x)=k\deg c_{n-1}(x)$. Thus, we have $c_{n-i}(x)\neq 0$ and $\deg c_{n-i}(x)=i\deg c_{n-1}(x)$ for all $1\leq i\leq n$.

It follows from equation \eqref{eq2.5} that $\deg c_1(x)=\deg d_0(x)+\deg c_0(x)$. It follows from equation \eqref{eq2.6} and Claim 2 that $\deg a_2(x)+2\deg c_{n-1}(x)=0$, which contradicts that $\deg a_2(x)\geq1$.

 If $n=1$, then it follows from equation \eqref{eq2.5} that $\deg c_1(x)=\deg d_0(x)+\deg c_0(x)$. Hence $\deg d_0(x)=0$, which contradicts equation \eqref{eq2.6}.
 Thus, $D$ has no Darboux polynomial. Whence $D$ is simple.
\end{proof}
\begin{prop}\label{prop2.2}
Let $D=y\partial _x+{(a_1(x)y+a_0(x))}{\partial _y}$ be a derivation of $K[x,y]$ with $a_1(x),a_0(x)\in K[x]$. If $D$ is simple, then $\operatorname{Im}D$ is not a Mathieu-Zhao space of $K[x,y]$.
\end{prop}
\begin{proof}
It follows from Theorem \ref{thm2.1} in \cite{15} that $a_0(x)\in K^*$, $\deg a_1(x)\geq1$. Let $a_0:=a_0(x)\in K^*$. Note that $1=a_0^{-1}D(y-\int a_1(x)dx)$. Then we have $1\in\operatorname{Im}D$. If $\operatorname{Im}D$ is a Mathieu-Zhao space, then $\operatorname{Im}D=K[x,y]$.

Claim : $x\notin \operatorname{Im}D$.

Suppose that $x\in \operatorname{Im}D$. Then there exists $G(x,y)\in K[x,y]$ such that
\begin{equation}
 yG_x+(a_1(x)y+a_0)G_y=x.\label{eq2.12}
\end{equation}

Let $G(x,y)=\sum_{i=0}^t{b_i(x)}y^i$ with $b_t(x)\neq0,b_i(x)\in K[x]$ for $0\leq i\leq t$. If $t=0$, then $D(G(x,y))=yb'_0(x)\neq x$. Thus, we have $t\geq1$. Then equation \eqref{eq2.12} has the following form :
\begin{equation}
y(\sum_{i=0}^t{b'_i(x)}y^i)+(a_1(x)y+a_0)(\sum_{i=1}^t{ib_i(x)}y^{i-1})=x.\label{eq4.13}
\end{equation}

We view the polynomials as in $K[x][y]$ with coefficients in $K[x]$ when we compare the coefficients of monomials on $y$. Comparing the coefficients of $y^{t+1}$ of equation \eqref{eq4.13}, we have
\begin{equation}
b'_t(x)=0.\label{eq4.15}
\end{equation}

Thus, we have $ b_t:=b_t(x)\in K^*$. Then we have the following equations:
\begin{equation}\label{eq4.14}
\left\{ \begin{array}{l   l}
b_{t-1}'(x)+ta_1(x)b_t=0 \\
b_{i-1}'(x)+ia_1(x)b_i(x)+(i+1)a_0b_{i+1}(x)=0\\
a_0b_1(x)=x
\end{array} \right.
\end{equation}
by comparing the coefficients of monomials on $y$ of equation \eqref{eq4.13} for $1\leq i\leq{t-1}$.

If $t\geq2$, then it follows from equations \eqref{eq4.14} that $\deg b_{t-i}(x)=i(\deg a_1(x)+1)$ for $1\leq i\leq t$. Thus, $\deg b_1(x)=(t-1)(\deg a_1(x)+1)\geq2$, which contradicts that $\deg b_1(x)=1$. If $t=1$, then it follows from the third equation of equations \eqref{eq4.14} that $\deg b_1(x) = 1$, which contradicts the fact that $b_1(x)\in K^*$. Hence we have $x\notin\operatorname{Im}D$. Therefore, $\operatorname{Im}D$ is not a Mathieu-Zhao space of $K[x,y]$.
\end{proof}
\begin{cor}\label{cor2.3}
Let $D=y\partial_x+(a_1(x)y+a_0)\partial_y$ be a derivation of $K[x,y]$ with $a_1(x)\in K[x],a_0\in K$. Then $D$ is simple iff $\operatorname{Im}D$ is not a Mathieu-Zhao space of $K[x,y]$.
\end{cor}
\begin{proof}
$``\Rightarrow"$ It follows from Proposition \ref{prop2.2}.

$``\Leftarrow"$
It suffices to prove that if $\operatorname{Im}D$ is not a Mathieu-Zhao space of $K[x,y]$, then $\deg a_1(x)\geq1,a_0\in K^*$. Next we prove that if $\deg a_1(x)\leq 0,a_0\in K^*$ or $a_0=0$, then $\operatorname{Im}D$ is a Mathieu-Zhao space. If $a_0=0$, let $D_1:=\partial_x+a_1(x)\partial_y$, then $\operatorname{Im}D_1=K[x,y]$. Hence $\operatorname{Im}D=(y)$. In particular, $\operatorname{Im}D$ is a Mathieu-Zhao space of $K[x,y]$.
If $\deg a_1(x)\leq 0,a_0\in K^*$, let $a_1:=a_1(x)\in K$, then $D=y\partial_x+(a_1y+a_0)\partial_y$ is locally finite. It follows from Theorem 3.1 in \cite{16} that $\operatorname{Im}D$ is a Mathieu-Zhao space of $K[x,y]$. Then the conclusion follows.
\end{proof}

\section{Some lemmas}

In this section, we prove some lemmas which we need in section 4.

\begin{lem}\label{lem3.1}
Let $D=y\partial _x+{(a_2(x)y^2+a_1(x)y+a_0)}{\partial _y}$ be a derivation of $K[x,y]$ with $a_2(x)\in K[x],\deg a_2(x)\geq1$ and $a_0\in K^*,a_1(x)\in K[x] \backslash \{0\}$. If $D$ has a Darboux polynomial $F(x,y)=\sum_{i=0}^n{c_i(x)}y^i\notin K$ with $c_n(x)\neq 0,c_i(x)\in K[x]$ for $0\leq i \leq n$ such that $D(F(x,y))=\Lambda(x,y)F(x,y)$ for some $\Lambda(x,y)\in K[x,y]$, then we have

$(1)$ $n\geq 1,c_n(x)\in K^*,\Lambda(x,y)=d_1(x)y+d_0(x)$ with $d_1(x)=na_2(x)$ and $d_0(x)\in K[x]$.

$(2)$
\begin{align}
c_{n-1}'(x)&=a_2(x)c_{n-1}(x)+(d_0(x)-na_1(x))c_n(x),\label{eq3.1}\\
(i+1)a_0c_{i+1}(x)&=(n-i+1)a_2(x)c_{i-1}(x)+(d_0(x)-ia_1(x))c_i(x)-c_{i-1}'(x),\label{eq3.2}\\
a_0c_1(x)&=d_0(x)c_0(x)\label{eq3.3}
\end{align}
for $1\leq i\leq{n-1}$ and $n\geq2$. If $n=1$, then we have equations \eqref{eq3.1}, \eqref{eq3.3}.

$(3)$ $c_{n-1}(x)\neq 0, d_0(x)-na_1(x)\notin K$ for $n\geq2$.

In particular, if $\deg a_1(x)=0$ and $n\geq 2$, then we have $\deg d_0(x)\geq1$ and $c_1(x)c_0(x)\neq 0$.
\end{lem}
\begin{proof}
If $n=0$, then $F(x,y)=c_0(x)$, $yc'_0(x)=\Lambda(x,y)c_0(x)$. Thus, we have $c_0(x)\in K$, which is a contradiction. Hence we have $n\geq 1$.

It follows from the equation $D(F(x,y))=\Lambda(x,y)F(x,y)$ that
\begin{equation}\label{eq3.4}
y{\sum_{i=0}^n c_i'(x) y^i}+(a_2(x)y^2+a_1(x)y+a_0){\sum_{i=1}^n ic_i(x) y^{i-1}}=\Lambda(x,y){\sum_{i=0}^n c_i(x) y^i}.
\end{equation}

Comparing the highest degree of $y$ of equation \eqref{eq3.4}, we have $\deg_y\Lambda(x,y)\leq 1$. Let $\Lambda(x,y)=d_1(x)y+d_0(x)$ with $d_1(x),d_0(x)\in K[x]$. Then equation \eqref{eq3.4} has the following form:
\begin{equation}
y{\sum_{i=0}^n c_i'(x) y^i}+(a_2(x)y^2+a_1(x)y+a_0){\sum_{i=1}^n ic_i(x) y^{i-1}}=(d_1(x)y+d_0(x)){\sum_{i=0}^n c_i(x) y^i}. \label{eq3.5}
\end{equation}

We view the polynomials as in $K[x][y]$ with coefficients in $K[x]$ when we comparing the coefficients of monomials on $y$. Comparing the coefficients of $y^{n+1}$ of equation \eqref{eq3.5}, we have
\begin{equation}
c'_n(x)+na_2(x)c_n(x)=d_1(x)c_n(x).\label{eq3.14}
\end{equation}

Thus, we  have $c_n(x)\in K^*$ and $d_1(x)=na_2(x)$ by comparing the degree of $x$ of the above equation. Comparing the  the coefficients of monomials on $y$ of equation \eqref{eq3.5}, we have the following equations:
\begin{align}
c_{n-1}'(x)&=a_2(x)c_{n-1}(x)+(d_0(x)-na_1(x))c_n(x),\label{eq3.6}\\
(i+1)a_0c_{i+1}(x)&=(n-i+1)a_2(x)c_{i-1}(x)+(d_0(x)-ia_1(x))c_i(x)-c_{i-1}'(x),\label{eq3.7}\\
a_0c_1(x)&=d_0(x)c_0(x)\label{eq3.8}
\end{align}
for $1\leq i\leq{n-1}$ and $n\geq2$. If $n=1$, then we have equations \eqref{eq3.6},\eqref{eq3.8}.

Suppose that $d_0(x)-na_1(x)=0$. Then it follows from equation \eqref{eq3.6} that $c_{n-1}(x)=0$. It follows from equation \eqref{eq3.7}$(i=n-1)$ that $na_0c_{n}(x)=-c_{n-2}'(x)\\+2a_2(x)c_{n-2}(x)$. If $c_{n-2}(x)=0$, then it follows from equation \eqref{eq3.7}$(i=n-1)$ that $c_{n}(x)=0$, which is a contradiction. If $c_{n-2}(x)\neq0$, then we have a contradiction by comparing the degree of $x$ of equation \eqref{eq3.7}$(i=n-1)$. Hence $d_0(x)-na_1(x)\neq0$. Suppose that $d_0(x)-na_1(x)\in K^*$. Then we have a contradiction by comparing the degree of $x$ of equation \eqref{eq3.6}. Thus, we have $d_0(x)-na_1(x)\notin K$. It follows from equation \eqref{eq3.6} that $c_{n-1}(x)\neq0$.

If $\deg a_1(x)=0$ and $c_1(x)=0$, then it follows from equation \eqref{eq3.8} that $c_0(x)=0$. It follows from equation \eqref{eq3.7} that $c_i(x)=0$ for $2\leq i\leq n$. This contradicts the fact that $c_n(x)\in K^*$. Analogously, we have $c_0(x)\neq 0$.

If $\deg a_1(x)=0$ and $ d_0(x)\in K$, then $d_0(x)-na_1(x)\in K$, a contradiction. Hence $\deg d_0(x)\geq1$.
\end{proof}

In order to prove Lemma \ref{lem3.3}, we first give the following lemma.

\begin{lem} \label{lem3.2}
If $c_i(x)\in K[x]$ satisfies the equations \eqref{eq3.2},\eqref{eq3.3} for $0\leq i\leq n$ and $n\geq2$, where $\deg a_1(x)>\deg a_2(x)\geq1$, $\deg c_n(x)=0,c_0(x)c_1(x)c_{n-1}(x)\neq0$ and $a_0\in K^*$, then $d_0(x)\neq ja_1(x)$ for all $1\leq j\leq n-1$.
\end{lem}
\begin{proof}
Suppose that $d_0(x)=a_1(x)$. If $n=2$, then we have a contradiction by comparing the degree of $x$ of equation \eqref{eq3.2}($i=1$). If $n\geq3$, then substituting equations \eqref{eq3.3}, \eqref{eq3.2}$(i=1)$ into equation \eqref{eq3.2}$(i=2)$, we have
\begin{equation}
3a_0c_3(x)=(\frac{n}{2}-1)\frac {a_1(x)a_2(x)c_0(x)}{\displaystyle  a_0}-\frac {a_1'(x)c_0(x)}{\displaystyle  a_0}-\frac {a_1(x)c_0'(x)}{\displaystyle  2a_0}. \label{eq3.9}
\end{equation}

It follows from equation \eqref{eq3.2}$(i=1)$ that $\deg c_2(x)=\deg a_2(x)+\deg c_0(x)$. Then it follows from equation \eqref{eq3.9} that $\deg c_3(x)=\deg a_1(x)+\deg c_2(x)>\deg c_2(x)\geq 1$. If $n\geq4$, then it follows from equation \eqref{eq3.2}$(i=n-1)$ that $\deg c_{n-2}(x)+\deg a_2(x)=\deg c_{n-1}(x)+\deg a_1(x)$. Then we have $\deg c_{n-2}(x)>\deg c_{n-1}(x)$. Continuing in this way we arrive at $\deg c_2(x)>\deg c_3(x)$, which is a contradiction. If $n=3$, then we have $\deg c_3(x)=0$, which is a contradiction.

Suppose that $d_0(x)= i_0a_1(x)$ with $2\leq i_0\leq n-1$. Then $\deg(d_0(x)-ia_1(x))=\deg a_1(x)>\deg a_2(x) $ for $1\leq i\leq n$ and $i\neq i_0$.
It follows from equation \eqref{eq3.2} that $\deg c_{i+1}(x)=\deg a_1(x)+\deg c_i(x)$ for $1\leq i\leq i_0-1$ and $\deg c_{i_0+1}(x)=\deg a_2(x)+\deg c_{i_0-1}(x)$. Then we have $\deg c_{i_0+1}(x)+\deg a_1(x)=\deg a_1(x)+\deg a_2(x)+\deg c_{i_0-1}(x)=\deg a_2(x)+\deg c_{i_0}(x)$.

Substituting equation \eqref{eq3.2}$(i=i_0-1)$ and equation \eqref{eq3.2}$(i=i_0)$ into equation \eqref{eq3.2}$(i=i_0+1)$, we have
\begin{align}
(i_0+2)a_0c_{i_0+2}(x)=\frac {n-i_0}{\displaystyle i_0 a_0} a_2(x) [(n-i_0+2)a_2(x)c_{i_0-2}(x)+ a_1(x)c_{i_0-1}(x)\\\notag-c_{i_0-2}'(x)]
+\frac {-a_1(x)}{\displaystyle (i_0+1) a_0}[(n-i_0+1)a_2(x)c_{i_0-1}(x)-c_{i_0-1}'(x)]-c_{i_0}'(x).
\end{align}

Clearly,
\begin{align}
\frac {n-i_0}{\displaystyle i_0 a_0} a_2(x) a_1(x)c_{i_0-1}(x)+\frac {-a_1(x)}{\displaystyle (i_0+1) a_0}(n-i_0+1)a_2(x)c_{i_0-1}(x)\\ \notag
=\frac {n-2i_0}{\displaystyle i_0(i_0+1) a_0}a_1(x)a_2(x)c_{i_0-1}(x). \notag
\end{align}

If $n=2i_0$, then we consider
\begin{equation}
\frac {n-i_0}{\displaystyle i_0 a_0} a_2(x) (n-i_0+2)a_2(x)c_{i_0-2}(x)+\frac {a_1(x)}{\displaystyle (i_0+1) a_0}c_{i_0-1}'(x)-c_{i_0}'(x).\label{eq3.11}
\end{equation}

Substituting equation \eqref{eq3.2}$(i=i_0-1)$ into equation \eqref{eq3.11}, we have
\begin{align}
\frac {i_0+2}{\displaystyle a_0} a_2^2(x)c_{i_0-2}(x)+\frac {a_1(x)}{\displaystyle (i_0+1) a_0}c_{i_0-1}'(x)-
\frac {1}{\displaystyle i_0a_0 }[(i_0+2)a_2'(x)c_{i_0-2}(x)\\\notag+(i_0+2)a_2(x)c_{i_0-2}'(x)+a_1'(x)c_{i_0-1}(x)+a_1(x)c_{i_0-1}'(x)-c_{i_0-2}''(x)]\\\notag
=\frac {i_0+2}{\displaystyle a_0} (a_2^2(x)-\frac {a_2'(x)}{\displaystyle i_0})c_{i_0-2}(x)-\frac {1}{\displaystyle i_0 a_0}(i_0+2)a_2(x)c_{i_0-2}'(x)+\frac {1}{\displaystyle i_0 a_0}c_{i_0-2}''(x)\\ \notag-\frac {a_1'(x)}{\displaystyle i_0 a_0}c_{i_0-1}(x)-\frac {a_1(x)}{\displaystyle i_0(i_0+1) a_0}c_{i_0-1}'(x).
\end{align}

Since $2\deg a_2(x)+\deg c_{i_0-2}(x)-\deg c_{i_0-1}(x)+1-\deg a_1(x)
=2\deg a_2(x)-2\deg a_1(x)+1<0$ and the coefficient of $x$ with highest degree in $\frac {a_1'(x)}{ i_0 a_0}c_{i_0-1}(x)+\frac {a_1(x)}{i_0(i_0+1) a_0}c_{i_0-1}'(x)$ is nonzero, we have $\deg c_{i_0+2}(x)=\deg c_{i_0-1}(x)-1+\deg a_1(x)
\geq\deg c_{i_0-1}(x)+\deg a_2(x)=\deg c_{i_0+1}(x)$.
 If $i_0\geq 3$, then it follows from equation \eqref{eq3.2} that $\deg c_{i+1}(x)=\deg c_i(x)+\deg a_1\geq1$ for $i_0+2\leq i\leq n-1$. Thus, we have $\deg c_{n}(x)=\deg c_{n-1}(x)+\deg a_1\geq1$, which is a contradiction. If $i_0=2$, then $n=4$, $\deg c_{4}(x)
\geq\deg c_{3}(x)\geq\deg a_2(x)+\deg c_1(x)\geq1$, which is a contradiction.

If $n\neq 2i_0$, then it follows from equation \eqref{eq3.2} that $\deg c_{i}(x)=\deg a_1(x)+\deg c_{i-1}(x)$ for $2\leq i\leq n$ and $i\neq i_0+1$, $\deg c_{i_0+1}(x)=\deg a_2(x)+\deg c_{i_0-1}(x)$. Thus, we have $\deg c_{n}(x)=\deg c_{n-1}(x)+\deg a_1\geq1$ or $\deg c_{n}(x)=\deg c_{n-2}(x)+\deg a_2\geq1$, which is a contradiction.
\end{proof}

\begin{lem} \label{lem3.3}
Let $D=y\partial _x+{(a_2(x)y^2+a_1(x)y+a_0)}{\partial _y}$ be a derivation of $K[x,y]$ with $a_1(x)\in K[x],a_2(x)\in K[x],\deg a_1(x)\geq1,\deg a_2(x)\geq1$ and $a_0\in K^*$. If $D$ has a Darboux polynomial $F(x,y)=\sum_{i=0}^n{c_i(x)}y^i\notin K$ with $c_n(x)\neq 0,c_i(x)\in K[x]$ for $0\leq i \leq n$ and $n\geq 2$ such that $D(F(x,y))=\Lambda(x,y)F(x,y)$ for some $\Lambda(x,y)\in K[x,y]$, then $d_0(x)\in K^*$ and $c_1(x)c_0(x)\neq 0$, where $d_0(x)$ be as in Lemma $\ref{lem3.1}$ $(1)$.
\end{lem}
\begin{proof}
Suppose that $d_0(x)=0$. Then it follows from equation \eqref{eq3.3} that $c_1(x)=0$.
If $n=2$, then we have a contradiction by comparing the degree of $x$ of equation \eqref{eq3.2}$(i=1)$.
If $n\geq3$ and $c_2(x)=0$, then it follows from equations \eqref{eq3.2} that $c_i(x)=0,3\leq i\leq n$, this contradicts the fact that $c_n(x)\in K^*$. If $n\geq3$ and $c_2(x)\neq 0$, then it follows from equation \eqref{eq3.1} that $c_{n-1}(x)\neq0$. Otherwise, $na_{1}(x)c_{n}(x)=0$, a contradiction. Thus, we have $\deg a_1(x)=\deg a_2(x)+\deg c_{n-1}(x)\geq\deg a_2(x)$. It follows from equations \eqref{eq3.2} that $ c_{i+1}(x)\neq0$ and $\deg c_{i+1}(x)=\deg a_1(x)+\deg c_i(x)$ for $2\leq i\leq n-1$. Then we have $\deg c_n(x)\geq 1$, which is a contradiction.

 Suppose that $c_1(x)=0$. Then it follows from equation \eqref{eq3.3} that $c_0(x)=0$. It follows from equations \eqref{eq3.2} that $c_i(x)=0$ for $2\leq i\leq n$. This contradicts the fact that $c_n(x)\in K^*$. Analogously, we have $c_0(x)\neq 0$.

 Let $t:=\max\{\deg d_0(x),\deg a_1(x)\}$. Suppose that $\deg d_0(x)\geq1$. It follows from $c_{n-1}(x)\neq0,d_0(x)-na_1(x)\notin K$(by Lemma \ref{lem3.1}) and equation \eqref{eq3.1} that
\begin{equation}
\deg (d_0(x)-na_1(x))=\deg a_2(x)+\deg c_{n-1}(x). \label{eq3.12}
\end{equation}

We distinguish four cases:

Case 1: If $\deg (d_0(x)-ia_1(x))=t$ for $1\leq i\leq n$, then $\deg (d_0(x)-na_1(x))=t$.
It follows from equation \eqref{eq3.12} that $t=\deg a_2(x)+\deg c_{n-1}(x)\geq\deg a_2(x)$. Thus, $\deg(d_0(x)-ia_1(x))\geq\deg a_2(x) $ for $1\leq i\leq n$.
It follows from equation \eqref{eq3.3} that $\deg c_1(x)=\deg c_0(x)+\deg d_0(x)> \deg c_0(x)$.
It follows from equations \eqref{eq3.2} that $\deg c_{i+1}(x)=\deg c_{i}(x)+\deg (d_0(x)-ia_1(x))=\deg c_{i}(x)+t$ for $1\leq i\leq n-1$.
Then we have $\deg c_n(x)=\deg c_{n-1}(x)+t\geq1$. This contradicts the fact that $c_n(x)\in K^*$.

Case 2: If there exists a $i_0\in \{1,...,n-1\}$ such that $d_0(x)-i_0a_1(x)=0$, then $\deg(d_0(x)-na_1(x))= \deg a_1(x)$.
It follows from equation \eqref{eq3.12} that $\deg a_1(x)=\deg a_2(x)+\deg c_{n-1}(x)\geq\deg a_2(x)$. Thus, $\deg(d_0(x)-ia_1(x))=\deg a_1(x)\geq\deg a_2(x) $ for $1\leq i\leq n$ and $i\neq i_0$.

It follows from Lemma \ref{lem3.2} that $\deg a_1(x)\leq\deg a_2(x)$. Thus, $\deg a_1(x)=\deg a_2(x)$. It follows from equation \eqref{eq3.12} that $\deg c_{n-1}(x)=0$. If $n=2$, then $\deg c_1(x)=0$. If $n\geq3$, then it follows from equation \eqref{eq3.2}$(i=n-1)$ that $i_0\neq n-1$. Otherwise, we have a contradiction by comparing the degree of $x$. Hence we have $\deg c_{n-2}(x)=\deg c_{n-1}(x)=0$ by comparing the degree of $x$ of the equation \eqref{eq3.2}$(i=n-1)$. Continuing in this way we arrive at $\deg c_1(x)=\deg c_2(x)=0$. It follows from equation \eqref{eq3.3} that $\deg c_1(x)=\deg d_0(x)+\deg c_0(x)\geq1$, which is a contradiction.

Case 3: If there exists a $i_0\in \{1,...,n-1\}$ such that $d_0(x)-i_0a_1(x)\in K^*$, then it follows from equation \eqref{eq3.12} that $\deg a_1(x)=\deg a_2(x)+\deg c_{n-1}(x)\geq\deg a_2(x)$. Thus, $\deg(d_0(x)-ia_1(x))=\deg a_1(x)\geq\deg a_2(x) $ for $1\leq i\leq n$ and $i\neq i_0$.

If $\deg a_1(x)>\deg a_2(x)$ and $n\geq3$, then it follows from equations \eqref{eq3.2} that $\deg c_{i_0+1}(x)=\deg c_{i_0}(x)$, $\deg c_{i+1}(x)=\deg c_i(x)+\deg a_1(x)$ for $i\neq i_0$ and $1\leq i\leq n-1$. Then we have $\deg c_{n}(x)\geq\deg c_1(x)+\deg a_1(x)\geq1$, which is a contradiction.

If $\deg a_1(x)>\deg a_2(x)$ and $n=2$, then it follows from equation \eqref{eq3.2} that $\deg c_2(x)=\deg c_1(x)$. It follows from equation \eqref{eq3.3} that $\deg c_1(x)=\deg c_0(x)+\deg d_0(x)\geq1$. Then we have $\deg c_2(x)\geq1$, which is a contradiction.

If $\deg a_1(x)=\deg a_2(x)$, then it follows from equation \eqref{eq3.1} that $\deg c_{n-1}(x)=0$. If $n=2$, then we have $\deg c_1(x)=0$. If $n\geq3$ and $i_0=n-1$, then it follows from equation \eqref{eq3.2}$(i=n-1)$ that $c_{n-2}(x)=0$. Otherwise, we have a contradiction by comparing the degree of $x$ of equation \eqref{eq3.2}$(i=n-1)$. If $c_{n-3}(x)=0$, then it follows from equation \eqref{eq3.2}$(i=n-2)$ that $(n-1)a_0c_{n-1}(x)=0$, a contradiction. If $c_{n-3}(x)\neq0$, then we have a contradiction by comparing the degree of $x$ of equation \eqref{eq3.2}$(i=n-2)$. Thus, we have $i_0\neq n-1$. It follows from equation \eqref{eq3.2}$(i=n-1)$ that $\deg c_{n-2}(x)=0$. Continuing in this way we arrive at $i_0\neq 2$ and $\deg c_1(x)=0$. It follows from equation \eqref{eq3.3} that $\deg c_1(x)\geq1$, which is a contradiction.

Case 4:  If $\deg (d_0(x)-i a_1(x))\geq1$ for $1\leq i\leq n$ and there exists a $i_0\in \{1,...,n\}$ such that $\deg (d_0(x)-i_0a_1(x))<t$, then $\deg a_1(x)=\deg d_0(x)$ and $\deg (d_0(x)-i a_1(x))=\deg a_1(x)=\deg d_0(x)$ for $1\leq i\leq n$ and $i\neq i_0$.

If $i_0=n$, then $\deg (d_0(x)-na_1(x))<\deg a_1(x)$. It follows from equation \eqref{eq3.12} that $\deg a_1(x)>\deg c_{n-1}(x)+\deg a_2(x)\geq\deg a_2(x)$, then $\deg (d_0(x)-ia_1(x))=\deg a_1(x)>\deg a_2(x) $ for $1\leq i\leq n-1$. It follows from equation \eqref{eq3.3} that $\deg c_1(x)\geq \deg c_0(x)+1$. It follows from equations \eqref{eq3.2} that $c_i(x)\neq0$ and $\deg c_{i+1}(x)=\deg c_i(x)+\deg a_1(x)$ for $1\leq i\leq n-1$. Then $\deg c_n(x)=\deg c_{n-1}(x)+\deg a_1(x)\geq1$, which is a contradiction.

If $i_0\neq n$, then $\deg (d_0(x)-na_1(x))=\deg a_1(x)$.
It follows from equation \eqref{eq3.12} that $\deg a_1(x)=\deg c_{n-1}(x)+\deg a_2(x)$. Then we have $\deg (d_0(x)-ia_1(x))=\deg a_1(x)\geq \deg a_2(x)$ for $i\neq i_0$ and $1\leq i\leq n$.

If $\deg a_1(x)=\deg a_2(x)$, then $\deg c_{n-1}(x)=0$. It follows from equation \eqref{eq3.2}$(i=n-1)$ that $i_0\neq n-1$. Comparing the degree of $x$ of the equation \eqref{eq3.2}$(i=n-1)$, we have $\deg c_{n-2}(x)=0$. Similarly, we have $\deg c_{i}(x)=0$ for $1\leq i\leq n-2$. Thus, we have $\deg c_{i}(x)=0$ for $1\leq i\leq n-1$. It follows from equation \eqref{eq3.3} that $\deg c_1(x)=\deg d_0(x)+\deg c_0(x)\geq1$, which is a contradiction.

If $\deg a_1(x)>\deg a_2(x)$, then it follows from equation \eqref{eq3.3} that $\deg c_1(x)\geq \deg c_0(x)+1$. It follows from equation \eqref{eq3.2} that $\deg c_{i_0+1}(x)=\deg c_{i_0}(x)+\deg {(d_0(x)-i_0a_1(x))}\geq 1+\deg c_{i_0}(x)$. Then we have $c_i(x)\neq0$ and $\deg c_{i+1}(x)=\deg c_i(x)+\deg a_1(x)$ for $i\neq i_0, 1\leq i\leq n-1$. Thus, we have $\deg c_{n}(x)\geq\deg c_{n-1}(x)+1$, which is a contradiction.
\end{proof}

\section{The main result}

In this section, we prove that $D=y\partial_x+(a_2(x)y^2+a_1(x)y+a_0(x))\partial_y$ with $a_0(x), a_1(x), a_2(x)\in K[x]$ is simple iff $(1)$ $a_0(x)\in K^*$,
$(2)$ $ \deg a_1(x)\geq 1$ or $\deg a_2(x)\geq 1$,
$(3)$ there exist no $l\in K^*$ such that $a_2(x)=la_1(x)-l^2a_0(x)$.

\begin{thm} \label{thm4.1}
Let $D=y\partial _x+{(a_2(x)y^2+a_1y+a_0(x))}{\partial _y}$ be a derivation of $K[x,y]$ with $a_1\in K,a_0(x),a_2(x)\in K[x]$. Then $D$ is simple if and only if $a_0(x)\in K^*$ and $\deg a_2(x)\geq1$.
\end{thm}
\begin{proof}
$``\Rightarrow"$ If $a_0(x)\notin K^*$, then $(y,a_0(x))$ is a $D$-stable ideal. Thus, $D$ is not simple. Hence we have $a_0(x)\in K^*$.
If $a_2(x)=0,a_1=0$, then $(\frac{1}{2}y^2-xa_0(x))$ is a $D$-stable ideal. Thus, $D$ is not simple.
If $a_2(x)=0,a_1\in K^*$, then $(y+\frac{a_0(x)}{a_1})$ is a $D$-stable ideal. Thus, $D$ is not simple.
If $a_2(x)\in K^*$, then $(y^2+\frac{a_1}{a_2(x)}y+\frac{a_0(x)}{a_2(x)})$ is a $D$-stable ideal. Thus, $D$ is not simple. Hence we have $\deg a_2(x)\geq1$.

$``\Leftarrow"$
Let $a_0:=a_0(x)\in K^*$. It follows from Proposition 2.1 in \cite{17} that it suffices to prove that $(y,a_2(x)y^2+a_1y+a_0)=(1)$ and $D$ has no Darboux polynomial. Clearly, $(y,a_2(x)y^2+a_1y+a_0)=(a_0)=(1)$. If $a_1=0$, then it follows from Theorem \ref{thm2.1} that $D$ is simple. Then we can assume $a_1\in K^*$.
Suppose that $F(x,y)\in K[x,y]$ is a Darboux polynomial of $D$ and $F(x,y)={\sum_{i=0}^n c_i(x) y^i}\notin K$ with $c_n(x)\neq 0,c_i(x)\in K[x]$ for $0\leq i \leq n$. Then we have $D(F(x,y))=\Lambda(x,y)F(x,y)$ for some $\Lambda(x,y)\in K[x,y]$.

It follows from Lemma \ref{lem3.1} that $n\geq1,\Lambda(x,y)=d_1(x)y+d_0(x)$ with $d_1(x),d_0(x)\\\in K[x],c_n(x)\in K^*$ and $d_1(x)=na_2(x)$.

If $n=1$, then it follows from Lemma \ref{lem3.1} that
\begin{equation}
\left\{ \begin{array}{l   l}
c_{0}'(x)=a_2(x)c_{0}(x)+(d_0(x)-a_1)c_1(x),\label{eq4.3}\\
a_0c_1(x)=d_0(x)c_0(x).
\end{array} \right.
\end{equation}

It follows from equations \eqref{eq4.3} that $\deg c_0(x)=\deg d_0(x)=0$, $d_0(x)=a_1, c_0(x)=0$ and $c_1(x)=0$, which is a contradiction.

 If $n\geq2$, then it follows from Lemma \ref{lem3.1} that we have the following equations:
\begin{align}
c_{n-1}'(x)&=a_2(x)c_{n-1}(x)+(d_0(x)-na_1)c_n(x),\label{eq4.4}\\
(i+1)a_0c_{i+1}(x)&=(n-i+1)a_2(x)c_{i-1}(x)+(d_0(x)-ia_1)c_i(x)-c_{i-1}'(x),\label{eq4.5}\\
a_0c_1(x)&=d_0(x)c_0(x)\label{eq4.6}
\end{align}
for $1\leq i\leq{n-1}$.\\
\indent It follows from Lemma \ref{lem3.1} that $\deg d_0(x)\geq1,c_0(x)c_1(x)\neq0$ and $c_{n-1}(x)\neq0$.\\
\indent It follows from equation \eqref{eq4.6} that $\deg c_1(x)=\deg d_0(x)+\deg c_0(x)\geq1+\deg c_0(x)$. It follows from equation \eqref{eq4.4} that $\deg (d_0(x)-na_1)=\deg c_{n-1}(x)+\deg a_2(x)$. Then we have $\deg d_0(x)=\deg (d_0(x)-ia_1)=\deg (d_0(x)-na_1)\geq\deg a_2(x)$ for $1\leq i\leq n$. It follows from equation \eqref{eq4.5} that $c_i(x)\neq0$ and $\deg c_{i+1}(x)=\deg c_i(x)+\deg d_0(x)$ for $1\leq i\leq n-1$. Thus, we have $\deg c_{n}(x)=\deg c_{n-1}(x)+\deg d_0(x)\geq1$, which is a contradiction. Hence $D$ has no Darboux polynomial. Therefore, $D$ is simple.
\end{proof}

\begin{thm} \label{thm4.2}
Let $D=y\partial _x+{(a_2(x)y^2+a_1(x)y+a_0(x))}{\partial _y}$ be a derivation of $K[x,y]$ with $a_0(x),a_1(x),a_2(x)\in K[x]$. Then $D$ is simple iff the following conditions hold:

$(1)$ $a_0(x)\in K^*$,

$(2)$ $\deg a_1(x)\geq1$ or $\deg a_2(x)\geq1$,

$(3)$ there exist no $l\in K^*$ such that $a_2(x)=la_1(x)-l^2a_0(x)$.
\end{thm}
\begin{proof}
$``\Rightarrow"$ If $a_0(x)\notin K^*$, then $(y,a_0(x))$ is a $D$-stable ideal. Thus, $D$ is not simple. Hence we have $a_0(x)\in K^*$.
If $\deg a_1(x)\leq0,\deg a_2(x)=0$, then $(y^2+\frac{a_1(x)} {a_2(x)}y+\frac {a_0(x)} { a_2(x)})$ is $D$-stable. Thus, $D$ is not simple.
If $\deg a_1(x)=0, a_2(x)=0$, then $(y+\frac {a_0(x)} { a_1(x)})$ is $D$-stable. Thus, $D$ is not simple.
If $a_1(x)=0,a_2(x)=0$, then $(\frac{1}{2}y^2-xa_0(x))$ is a $D$-stable ideal. Thus, $D$ is not simple.
Hence we have $\deg a_2(x)\geq1$ or $\deg a_1(x)\geq 1$.

Clearly, $(y,a_2(x)y^2+a_1(x)y+a_0(x))=(a_0(x))=(1)$. It follows from Proposition 2.1 in \cite{17} that $D$ is simple if and only if $D$ has no Darboux polynomial. If $a_2(x)=la_1(x)-l^2a_0(x)$ for some $l\in K^*$, then $D(y+l^{-1})=l\cdot (a_1(x)y-la_0(x)y+a_0(x))(y+l^{-1})$. Hence $(y+l^{-1})$ is $D$-stable, whence $D$ is not simple. Then the conclusion follows.

$``\Leftarrow"$
Let $a_0:=a_0(x)\in K^*$. If $a_2(x)\in K$, then it follows from Theorem 2.2 in \cite{15} that $D$ is simple. If $\deg a_2(x)\geq 1,a_1\in K$, then it follows from Theorem \ref{thm4.1} that $D$ is simple. Thus, we can assume $\deg a_1(x)\geq 1$ and $\deg a_2(x)\geq 1$ in the following arguments.

 It follows from Proposition 2.1 in \cite{17} that it suffices to prove that $(y,a_2(x)y^2+a_1(x)y+a_0)=(1)$ and $D$ has no Darboux polynomial. Clearly, $(y,a_2(x)y^2+a_1(x)y+a_0)=(a_0)=(1)$.
Suppose that $F(x,y)\in K[x,y]$ is a Darboux polynomial of $D$ and $F(x,y)={\sum_{i=0}^n c_i(x) y^i}\notin K$ with $c_n(x)\neq 0,c_i(x)\in K[x]$ for $0\leq i \leq n$. Then we have $D(F(x,y))=\Lambda(x,y)F(x,y)$ for some $\Lambda(x,y)\in K[x,y]$.

It follows from Lemma \ref{lem3.1} that $n\geq1,\Lambda(x,y)=d_1(x)y+d_0(x)$ with $d_1(x),d_0(x)\\\in K[x],c_n(x)\in K^*$ and $d_1(x)=na_2(x)$.

If $n=1$, then we have the following equations:
\begin{equation}
\left\{ \begin{array}{l   l}
c_{0}'(x)=a_2(x)c_{0}(x)+(d_0(x)-a_1(x))c_1(x),\label{eq4.8}\\
a_0c_1(x)=d_0(x)c_0(x)
\end{array} \right.
\end{equation}

It follows from the above equations that $0=\deg c_{1}(x)=\deg c_0(x)+\deg d_0(x)$. Then $\deg c_0(x)=\deg d_0(x)=0$  and $c_0(x)\neq 0$. Hence we have $a_2(x)=\frac{c_1(x)}{c_0(x)}a_1(x)-\frac{c_{1}^2(x)}{c_0^2(x)}a_0(x)$, which is a contradiction.

If $n\geq 2$, then it follows from Lemma \ref{lem3.1} that
\begin{align}
c_{n-1}'(x)&=a_2(x)c_{n-1}(x)+(d_0(x)-na_1(x))c_n(x),\label{eq4.9}\\
(i+1)a_0(x)c_{i+1}(x)&=(n-i+1)a_2(x)c_{i-1}(x)+(d_0(x)-ia_1(x))c_i(x)-c_{i-1}'(x),\label{eq4.10}\\
a_0c_1(x)&=d_0(x)c_0(x)\label{eq4.11}
\end{align}
for $1\leq i\leq{n-1}$.

It follows from Lemma \ref{lem3.3} that $d_0(x)\in K^*$.

It follows from equation \eqref{eq4.9} that $c_{n-1}(x)\neq0$ and $\deg c_{n-1}(x)+\deg a_2(x)=\deg a_1(x)$. Then we have $\deg a_1(x)\geq\deg a_2(x)$.

If $\deg a_1(x)>\deg a_2(x)$, then it follows from equation \eqref{eq4.11} that $\deg c_1(x)=\deg c_0(x)$. It follows from equation \eqref{eq4.10} that $c_i(x)\neq0$ and $\deg c_{i+1}(x)=\deg c_i(x)+\deg a_1(x)\geq \deg a_1(x)\geq1$ for $1\leq i\leq n-1$. Then we have $\deg c_{n}(x)=\deg c_{n-1}(x)+\deg a_1(x)\geq1$. This contradicts the fact that $c_n(x)\in K^*$.

If $\deg a_1(x)=\deg a_2(x)$, then $\deg c_{n-1}(x)=0$. Comparing the degree of $x$ the equation \eqref{eq4.10}$(i=n-1)$, we have $\deg c_{n-2}(x)=0$. Similarly, we have $\deg c_{i}(x)=0$ for $0\leq i\leq n-2$. Thus, we have $\deg c_{i}(x)=0$ for $0\leq i\leq n$. Since $d_0(x)\in K^*$, it follows from equation \eqref{eq4.9} that $c_{n-1}(x)\neq 0$. Let $c_i:=c_i(x), d_0:=d_0(x)$ for $1\leq i\leq n$. It follows from equation \eqref{eq4.9} that $a_2(x)=\frac{nc_n}{c_{n-1}}a_1(x)-\frac{c_n}{c_{n-1}}d_0$. Then it follows from equation \eqref{eq4.10}$(i=n-1)$ that
\begin{equation}\label{eq4.12}
(\frac{2nc_nc_{n-2}}{c_{n-1}}-(n-1)c_{n-1})a_1(x)+d_0c_{n-1}-\frac{2c_nc_{n-2}}{c_{n-1}}d_0-na_0c_n=0.
\end{equation}
Then we have $c_{n-2}=\frac{(n-1)c_{n-1}^2}{2nc_n}$ and $d_0c_{n-1}-\frac{2c_nc_{n-2}}{c_{n-1}}d_0-na_0c_n=0$ by comparing the coefficients of $x$ with highest degree of equation \eqref{eq4.12}. Then we have $d_0=\frac{n^2c_n}{c_{n-1}}a_0$. Thus, we have $a_2(x)=\frac{nc_n}{c_{n-1}}a_1(x)-(\frac{nc_n}{c_{n-1}})^2a_0$, which is a contradiction. Hence $D$ has no Darboux polynomial. Therefore, $D$ is simple.
\end{proof}

\section{ The images of some derivations}

In this section, we prove that the images of some derivations of $K[x,y_1,\dots,y_n]$ and $K[y_1,\dots,y_n]$ are Mathieu-Zhao spaces.

\begin{thm}\label{thm5.1} Let $D=\partial _x+{\sum_{i=1}^n \gamma_i(x) y_i^{k_i}}{\partial _i}$ be a derivation of $K[x,y_1,\dots,y_n]$ and $k_i\geq 1,\gamma_i(x)\in K[x] \backslash \{0\}$ for all $1\leq i\leq n$. Then $\operatorname{Im}D$ is a Mathieu-Zhao space if and only if $k_i=1$ and $\gamma_i(x)\in K$ for all $1\leq i\leq n$.
\end{thm}
\begin{proof}
$``\Leftarrow"$  If $k_i=1$ and $\gamma_i(x)\in K$ for all $1\leq i\leq n$, then D is locally finite. Note that $1\in \operatorname{Im}D$ . Thus, the conclusion follows from Proposition 1.4 in \cite{10}.

$``\Rightarrow"$ Since $1\in \operatorname{Im}D$, we have $\operatorname{Im}D=K[x,y_1,\dots,y_n]$ if $\operatorname{Im}D$ is a Mathieu-Zhao space of $K[x,y_1,\dots,y_n]$.

If there exists $i_0\in \{1,2,\dots,n\}$ such that $k_{i_0}>1$, then we claim $y_{i_0}\notin \operatorname{Im}D$. Without loss of generality, we can assume that ${i_0}=1$. Suppose that ${y_1}\in \operatorname{Im}D$. Then there exists $f(x,y_1,\dots,y_n)\in K[x,y_1,\dots,y_n]$ such that
\begin{equation}
D(f(x,y_1,\dots,y_n))=y_1. \label{eq5.1}
\end{equation}

Let $f={\sum_{i=0}^t d_i(x) y_1^{i}+\sum_{|\alpha|=1}^d{c_\alpha(x)}y_1^{\alpha_1}\cdots y_n^{\alpha_n}}$ with $d_t(x)\neq 0$ and $\alpha_2+\alpha_3+\dots+\alpha_n\geq 1$, where ${c_\alpha(x)}\in K[x]$, $\alpha \in \mathbb{N}^n$,$|\alpha|=\alpha_1+\alpha_2+\alpha_3+\dots+\alpha_n$ and $d_i(x)\in K[x]$ for all $0\leq i\leq t$. It follows from equation \eqref{eq5.1} that
\begin{equation}
f_x+{\sum_{i=1}^n \gamma_i(x)y_i^{k_i}f_{y_i}}=y_1.\label{eq5.2}
\end{equation}

If $t\geq1$, then we have
\begin{equation}
{\sum_{i=0}^t d_i'(x) y_1^{i}}+\gamma_1(x)y_1^{k_1}\sum_{i=1}^t i d_i(x) y_1^{i-1}=y_1\label{eq5.3}
\end{equation}
by comparing the part of degree zero with respect to $y_2,\dots,y_n$ of equation \eqref{eq5.2}. Thus, we have $t \gamma_1(x)d_t(x)=0$ by comparing the coefficients of  $y_1^{k_1+t-1}$ of equation \eqref{eq5.3}, which is a contradiction.

If $t=0$, then we have $ d_0'(x)=y_1$
by comparing the part of degree zero with respect to $y_2,\dots,y_n$ of equation \eqref{eq5.2}, which is a contradiction.

If $k_i=1$ for all $1\leq i\leq n$ and there exists $j_0\in \{1,2,\dots,n\}$ such that $\deg{\gamma_{j_0}(x)}\geq 1$, then without loss of generality, we can assume that $j_0=1$, we claim $y_1\notin \operatorname{Im}D$. Suppose that $y_1\in \operatorname{Im}D$. Then there exists $ g(x,y_1,\dots,y_n)\in K[x,y_1,\dots,y_n]$ such that
\begin{equation}
D(g(x,y_1,\dots,y_n))=y_1.\label{eq5.4}
\end{equation}

Let $g=g^{(d)}+g^{(d-1)}+\dots+g^{(1)}+g^{(0)}$ with $g^{(d)}\neq0$, where $g^{(j)}$ is the homogeneous part of degree $j$ with respect to $y_1,\dots,y_n$ of $g$. It follows from equation (4.4) that
\begin{equation}
g_x+{\sum_{i=1}^n \gamma_i(x){y_i}}g_{y_i}=y_1. \label{eq5.5}
\end{equation}

Then we have
\begin{equation}
g_x^{(1)}+{\sum_{i=1}^n \gamma_i(x){y_i}}g_{y_i}^{(1)}=y_1. \label{eq5.6}
\end{equation}
by comparing the homogeneous part of degree one with respect to $y_1,\dots,y_n$ of equation \eqref{eq5.5}.

We view the polynomials as in $K[x][y_1,\dots,y_n]$ with coefficients in $K[x]$ when we comparing the coefficients of $y_1^{l_1}\dots y_n^{l_n}$. Let ${g^{(1)}}={c_1(x)y_1}+\dots+c_n(x)y_n$ with $c_i(x)\in K[x]$ for $1\leq i\leq n$. Then equation \eqref{eq5.6} has the following form :
\begin{equation}
\sum_{i=1}^n {{c_i}'(x)}{y_i}+{\sum_{i=1}^n {\gamma_i(x)}{y_i}{c_i}(x)}=y_1. \label{eq5.7}
\end{equation}

Thus, we have
\begin{equation}
{{c_1}'(x)}+\gamma_1(x)c_1(x)=1   \label{eq5.8}
\end{equation}
by comparing the coefficients of $y_1$ of equation \eqref{eq5.7}. Since $\deg \gamma_1(x)\geq 1$, we have $c_1(x)=0$ by comparing the degree of $x$ of equation \eqref{eq5.8}. Then equation \eqref{eq5.8} is 0=1, which is a contradiction. Then the conclusion follows.
\end{proof}

\begin{cor}\label{cor5.2}
Let $D={\partial_x}+{\sum_{i=1}^n \gamma_i(x) y_i^{k_i}}{\partial_i}$ be a derivation of $K[x,y_1,\dots,y_n]$ and $k_i\geq 1,\gamma_i(x)\in K[x]$ for all $1\leq i\leq n$. Then $\operatorname{Im}D$ is a Mathieu-Zhao space if and only if $D$ is locally finite.
\end{cor}
\begin{proof}
If $\gamma_i(x)=0$ for all $1\leq i\leq n$, then $D$ is locally finite and $\operatorname{Im}D$ is a Mathieu-Zhao space of $K[x,y_1,\dots,y_n]$.
If there exists $i_0\in \{1,2,\dots,n\}$ such that $\gamma_{i_0}(x)\neq0$, then the conclusion follows from Theorem \ref{thm5.1}.
\end{proof}

\begin{thm}\label{thm5.3}
Let $D={\sum_{i=1}^n \gamma_i y_i^{k_i}}{\partial _i}$ be a derivation of $K[y_1,\dots,y_n]$, $n\geq2$ and $\gamma_i\in K^*$ for all $1\leq i\leq n$.Then $\operatorname{Im}D$ is a Mathieu-Zhao space if and only if $k_i\leq1$ for all $1\leq i\leq n$.
\end{thm}
\begin{proof}
$``\Leftarrow"$ If $k_i=1$ for all $1\leq i\leq n$, then $D={\sum_{i=1}^n \gamma_i{y_i}}{\partial _i}$. Then the conclusion follows from Theorem 2.5 in \cite{9}.

If there exists $k_i=0$, then $1\in \operatorname{Im}D$ and $D$ is locally finite. Then the conclusion follows from Proposition 1.4 in \cite{10}.

$``\Rightarrow"$
Suppose that there exists $k_{i_0}>1$ and $k_{i}\neq0$ for all $1\leq i\leq n$, then we claim $y_{i_0}{y_{2}}^m\notin$ $\operatorname{Im}D$, where $m\gg 0$. Without loss of generality, we can assume that ${i_0}=1$. Suppose that $y_1 y_{2}^m\in$ $\operatorname{Im}D$. Then there exists $f(y_1,\dots,y_n)\in K[y_1,\dots,y_n]$ such that
\begin{equation}
D(f(y_1,\dots,y_n))=y_1 y_{2}^m.  \label{eq5.9}
\end{equation}

Let $f=h(y_1,y_2)+\sum_{|\alpha|=1}^p{c_\alpha}{y_1^{\alpha_1}\cdots y_n^{\alpha_n}}\neq0$ with $\alpha_3+\dots+\alpha_n\geq 1$, where ${c_\alpha}\in K$, $|\alpha|=\alpha_1+\alpha_2+\alpha_3+\dots+\alpha_n$ and $h(y_1,y_2)\in K[y_1,y_2]$.

It follows from equation \eqref{eq5.9} that
\begin{equation}
{\sum_{i=1}^n \gamma_i{y_i^{k_i}}f_{y_i}}=y_1 y_{2}^m.\label{eq5.10}
\end{equation}

Let $h(y_1,y_2)=\sum_{i=0}^t{c_i(y_2)}{y_1^{i}}$ with ${c_t(y_2)}\neq0$, where ${c_i(y_2)}\in K[y_2]$ for $0\leq i\leq t$. If $t=0$, then $D(h)\in K[y_2]$ and $D(h)\neq y_1y_{2}^m$. Thus we have $t\geq1$.

Then we have
\begin{equation}
\gamma_1{y_1^{k_1}}{\sum_{i=1}^t i{c_i(y_2)}{y_1^{i-1}}}+\gamma_2{y_2^{k_2}}{\sum_{i=0}^t {{c_i}'(y_2)}{y_1^{i}}}=y_1 y_{2}^{m} \label{eq5.11}
\end{equation}
by comparing the homogeneous part of degree zero with respect to $y_3,\dots,y_n$ of equation \eqref{eq5.10}. Thus, we have $t \gamma_1 c_t(y_2)=0$ by comparing the coefficients of  $y_1^{k_1+t-1}$ of equation \eqref{eq5.11}, which is a contradiction.

Note that $y_{2}^m=D(\frac {y_{2}^{m-k_2+1}} {\displaystyle {(m-k_2+1)a_2}}) $ for $m\geq k_2$, then we have $y_{2}^m\in$ $\operatorname{Im}D$ for $m\gg 0$, but $y_{1} y_{2}^m\notin$ $\operatorname{Im}D$ for $m\gg 0$, which is a contradiction.

Now suppose that there exists $k_{j_0}>1$ and $k_{p_0}=0$. Without loss of generality, we can assume that ${j_0}=1$ and ${p_0}=2$.

If $n=2$, then $D= \gamma_1 y_1^{k_1}{\partial _1}+ \gamma_2 {\partial _2}$. It follows from Theorem \ref{thm5.1} that $\operatorname{Im}D$ is not a Mathieu-Zhao space, which is a contradiction.

If $n\geq 3$, then $D=\gamma_1 y_1^{k_1}{\partial _1}+ \gamma_2 {\partial _2}+{\sum_{i=3}^n \gamma_i y_i^{k_i}}{\partial _i}$. Clearly, $1\in\operatorname{Im}D$. We have $\operatorname{Im}D  =K[y_1,\dots,y_n]$ if $\operatorname{Im}D$ is a Mathieu-Zhao space of $K[y_1,\dots,y_n]$.
We claim $y_1\notin \operatorname{Im}D$. Suppose that $y_1\in \operatorname{Im}D$. Then there exists $ g(y_1,\dots,y_n)\in K[y_1,\dots,y_n]$ such that
\begin{equation}
D(g(y_1,\dots,y_n))=y_1.\label{eq5.12}
\end{equation}

Let $g={\sum_{j=0}^s c_j y_1^{j}}$ with $c_{s}\neq0$, where $c_j\in K[y_2,\dots,y_n]$ for $0\leq j\leq s$. If $s=0$, then $D(g)\in K[y_2,\ldots,y_n]$ and $D(g)\neq y_1$. Thus we have $s\geq1$.

It follows from equation \eqref{eq5.12} that
\begin{equation}
\gamma_1 y_1^{k_1}{\sum_{j=1}^s jc_j y_1^{j-1}}
+\gamma_2 {\sum_{j=0}^s c_{j,2} y_1^{j}}
+{\sum_{i=3}^n \gamma_i y_i^{k_i}}{\sum_{j=0}^s c_{j,i} y_1^{j}}=y_1\label{eq5.13}
\end{equation}
where $c_{j,i}$ denotes $\frac{\partial c_j}{\partial y_i}$ for $2\leq i\leq n$ and $0\leq j\leq s $.

Thus, we have $s \gamma_1 c_s=0$ by comparing the coefficients of  $y_1^{k_1+s-1}$ of equation \eqref{eq5.13}, which is a contradiction. Then the conclusion follows. Similarly, we can prove the case that $k_{j_0}>1,k_{p_0}=\dots=k_{p_r}=0$.
\end{proof}

\section{ A conjecture}

Base on the conclusions of section 2 to section 4, we give the following conjecture.

\begin{conj}\label{conj6.1}
Let $D=y^{\alpha}\partial_x+(a_2(x)y^{\alpha+1}+a_1(x)y^{\alpha}+a_0(x))\partial_y$ be a derivation of $K[x,y]$ with $\alpha \in \mathbb{N}^*$. Then $D$ is simple if and only if

$(1)$ $a_0(x)\in K^*$,

$(2)$ $\deg a_2(x)\geq 1$ or $\deg a_1(x)\geq 1$,

$(3)$ there is no $l\in K^*$ such that $a_2(x)=la_1(x)+(-1)^{\alpha}l^{\alpha+1}a_0(x)$.
\end{conj}

\begin{rem}
In our paper, we have proved Conjecture \ref{conj6.1} if $\alpha=1$. In addition, one direction of Conjecture \ref{conj6.1} is easy to prove. Thus, we have the following proposition.
\end{rem}

\begin{prop}
Let $D=y^{\alpha}\partial_x+(a_2(x)y^{\beta+1}+a_1(x)y^{\beta}+a_0(x))\partial_y$ be a derivation of $K[x,y]$ with $\alpha, \beta\in \mathbb{N}^*$ and $\alpha\leq \beta$. If $D$ is simple, then we have the following statements:

$(1)$ $a_0(x)\in K^*$,

$(2)$ $\deg a_2(x)\geq 1$ or $\deg a_1(x)\geq 1$,

$(3)$ there is no $l\in K^*$ such that $a_2(x)=la_1(x)+(-1)^{\beta}l^{\beta+1}a_0(x)$.
\end{prop}
\begin{proof}
$(1)$ If $a_0(x)\notin K^*$, then the ideal $(y,a_0(x))$ is $D$-stable. Hence $D$ is not simple, which is a contradiction. Thus, we have $a_0(x)\in K^*$.

$(2)$ If $a_2(x)\in K^*$ or $a_1(x)\in K^*$, then the ideal $(a_2(x)y^{\beta+1}+a_1(x)y^{\beta}+a_0(x))$ is $D$-stable. If $a_2(x)=a_1(x)=0$, then the ideal $(\frac{1}{\alpha+1}y^{\alpha+1}-a_0\cdot x)$ is $D$-stable. Hence $D$ is not simple, which is a contradiction. Thus, we have $\deg a_2(x)\geq 1$ or $\deg a_1(x)\geq 1$.

$(3)$ If there is $l_0\in K^*$ such that $a_2(x)=l_0a_1(x)+(-1)^{\beta}l_0^{\beta+1}a_0(x)$, then the ideal $(y+l_0^{-1})$ is $D$-stable. Hence $D$ is not simple, which is a contradiction. Then the conclusion follows.
\end{proof}

\end{document}